\documentclass{article}
\usepackage{amssymb,amsmath,theorem,chicago,epsf}
\usepackage[]{titlesec}

\titleformat{\section}[hang]
  {\scshape\normalsize\filcenter}{\S\thesection}{1em}{}
\titleformat{\subsection}[hang]
  {\scshape\normalsize}{\S\thesubsection}{1em}{}
\titleformat{\subsubsection}[runin]
  {\itshape\normalsize}{\S\thesubsubsection}{.5em}{}[.]
\titlespacing{\subsubsection}{0pc}{*1}{1em}

{\theorembodyfont{\slshape}
\newtheorem{theorem}{\slshape Theorem}
\newtheorem{proposition}{\slshape Proposition}

\newcommand{\qedsymbol}{\rule{.3\baselineskip}{.35\baselineskip}}
\newenvironment{proof}{\begin{trivlist}\item[]{\it Proof.}}{\hspace*{\fill}$\qedsymbol$\end{trivlist}}

\newcommand{\erf}[1]{\expandafter(\ref{#1})} 
\newcommand{\thrf}[1]{\expandafter\ref{#1}}  
\newcommand{\ct}[1]{\expandafter\citeANP{#1}~[\expandafter\citeyearNP{#1}]}
\newcommand{\lb}[1]{\expandafter\label{#1}}
\newcommand{\elb}[1]{\expandafter\label{#1}}

\newcommand\mfk\mathfrak
\newcommand\mcl\mathcal
\newcommand\mbb\mathbb
\newcommand\mtl\mathit
\newcommand\mbf\mathbf

\newcommand\onm\operatorname
\newcommand{\set}[2]{\left\{\vphantom{\bigl(}\mskip1mu #1:#2\mskip1mu\right\}}
\newcommand{\sset}[1]{\left\{\vphantom{\bigl(}\mskip1mu#1\mskip1mu\right\}}

\title{Stability  by KAM confinement of certain wild,~nongeneric 
relative equilibria 
of underwater vehicles with coincident centers of~mass~and~bouyancy}
\author{George W. Patrick}
\date{April 2002}

\begin{document}\maketitle
\begin{abstract}
Purely rotational relative equilibria of an ellipsoidal underwater
vehicle occur at nongeneric momentum where the symplectic reduced
spaces change dimension.  The stability these relative equilibria
under momentum changing perturbations is not accessible by Lyapunov
functions obtained from energy and momentum.  A blow-up construction
transforms the stability problem to the analysis symmetry-breaking
perturbations of Hamiltonian relative equilibria.  As such, the
stability follows by KAM theory rather than energy-momentum
confinement.
\end{abstract}
\section*{Introduction}
The phase space $T\mtl{SE}(3)$ with Lagrangian
\begin{equation}\elb{00}
L(A,a,\Omega,v)\equiv\frac12\Omega^t\mbf I\Omega+\frac12v^t\mbf Mv
\end{equation}
approximately models the motion of a neutrally buoyant vehicle
submerged in an inviscid irrotational fluid (see~\ct{LeonardNE-1997.1}
and the references therein), in the case of coincident centers of mass
and buoyancy.  Here tangent vectors of $\mtl{SE}(3)$ are represented
by left translation and elements of $\mtl{SE}(3)$ parameterize the
configurations of the vehicle by embedding a reference vehicle into
the fluid.  $\mbf I$ and $\mbf M$ are constant, positive definite,
$3\times 3$ matrices that can be calculated from the shape and mass
distribution of the vehicle.  This system admits the $\mtl{SE}(3)$
symmetry of the left action of $\mtl{SE}(3)$ on itself.

For an ellipsoidal vehicle with principle axes of inertia along the
axes of symmetry of the ellipsoid, $\mbf I=\onm{diag}(I_1,I_2,I_3)$
and $\mbf M=\onm{diag}(M_1,M_2,M_3)$ (i.e.  $\mbf I$ and $\mbf M$ are
diagonal). If $M_1=M_2$ and $I_1=I_2$ (or similarly if $M_1=M_3$ and
$I_1=I_3$ etc.) then there is a further material symmetry of the
system: $\mtl{SO}(2)=\sset{\onm{exp}(\mbf k^\wedge\theta)}$ acts as a
subgroup of $\mtl{SE}(3)$ by inverse multiplication on the right.  In
the case that $\mbf I$ and $\mbf M$ are both constants of the identity
then the material symmetry is $\mtl{SO}(3)$.

Lie-Poisson reduction yields the Poisson
phase space $\mfk{se}(3)^*=\bigl\{(\pi,p)\bigr\}$ where $\pi=\mbf
I\Omega$ and $p=\mbf Mv$. The equations of motion are
\begin{equation}\elb{60}
\frac{d\pi}{dt}=\pi\times\Omega+p\times v,\qquad\frac{dp}{dt}=p\times\Omega,
\end{equation}
and by direct substitution, for each $\alpha_e\in\mbb R$,
\begin{equation*}
p_e^{\alpha_e}:\quad\pi=\alpha_e\mbf{k},\qquad p=0,
\end{equation*}
is an equilibrium of the Poisson reduced systems and hence a relative
equilibrium of the original system. The generator is 
\begin{equation*}
\Omega_e^{\alpha_e}\equiv\frac{\alpha_e}{I_3}\mbf k,
  \quad v_e^{\alpha_e}\equiv0,
\end{equation*}
so the relative equilibrium corresponds to a stationary vehicle
rotating about an axis of symmetry which is aligned with the vertical.
This article is concerned with the stability of these relative
equilibria in the case where $I_3$ is not an intermediate axis of
symmetry, i.e. assuming $I_1<I_2<I_3$ or $I_3<I_2<I_1$.

The symplectic leaves of $\mtl{se}(3)^*$ are as follows. On the
complement of $p=0$ lie the generic symplectic leaves, all
diffeomorphic to $TS^2$, and which are the level sets of the two
Casimirs $|p|$ and $\pi\cdot p$. Nongeneric leaves occur within the
set $p=0$ and are the level sets of the subcasimir $|\pi|$.  Thus the
relative equilibria $p_e^{\alpha_e}$ correspond to Lyapunov stable
equilibria on the (nongeneric) symplectic leaves of $\mtl{se}(3)^*$
since the energy has a definite critical point when restricted to
those leaves. Were the symmetry group to be compact this
\emph{leafwise stability} would imply stability of the equilibrium
modulo the isotropy group of the momentum.  $\mtl{SE}(3)$, of course,
is not compact. The question is whether or not $p_e^{\alpha_e}$ are
stable under perturbations from nongeneric leaves into nearby generic
leaves.

\ct{LeonardNEMarsdenJE-1997.1} have identified this question as
particularly delicate, and the theory
of~\ct{PatrickGWRobertsRMWulffC-2001.1}, the sharpest possible for the
problem of establishing the stability of relative equilibria by
energy-momentum confinement in the case of noncompact symmetry,
corroborates that opinion. Patrick \emph{et al.}  separate generators
of relative equilibria into two complementary classes, \emph{tame} and
\emph{wild}.  The generator of $p_e^{\alpha_e}$ is tame if and only if
$\alpha_e=0$, corresponding to a stationary, nonrotating vehicle, in
which case $\mtl{SE}(3)$-stability follows directly, since the energy
has zero derivative at $p_e^{\alpha_e}$ and has positive definite
Hessian there.  However, if $\alpha_e\ne0$, the generator is wild and
the theory does not imply stability.

So it is an open question whether the relative equilibria
$p_e^{\alpha_e}$, $\alpha_e\ne0$, are $\mtl{SE}(3)$-stable or not, and
the problem appears inaccessible by energy-momentum confinement. This is
due to the presence of a noncompact symmetry group and
wild generators. 

\section{The blow-up construction}
In its essence the stability issue is one of perturbation from a
nongeneric symplectic leaf to nearby, higher dimensional, generic
leaves.  In order to bridge to Hamiltonian perturbation theory, which
is usually cast in a setting of a fixed canonical phase space, it is
natural to begin by normalizing the generic leaves.  The leaf
corresponding to
\begin{equation*}
|p|=a,\quad\pi\cdot p=b,
\end{equation*}
for $a>0$ is diffeomorphic to 
\begin{equation*}
TS^2=\set{(w,\dot w)\in\mbb R^3\times\mbb R^3}{|w|=1,w\cdot\dot w=0}
\end{equation*}
by the map
\begin{equation*}
w=\frac p{|p|},\quad\dot w=\pi-\frac{\pi\cdot p}{|p|^2}p,
\end{equation*}
the inverse map being given by, for fixed $a>0$ and $b\in\mbb R$,
\begin{equation}\elb{4}
p=aw,\quad\pi=\dot w+\frac ba w.
\end{equation}
Having normalized the symplectic leaves to the constant manifold
$TS^2$, one seeks to extend this to the nongeneric leaves within
$p=0$, which means extending it to $a=0$, since $p=aw$.  As it
stands~\erf{4} is poorly defined for $a=0$, but for fixed ratios of
$b/a$ it is well defined even for arbitrarily small $a$, suggesting
that the proper way to approach the nongeneric leaves from generic
ones is through constant $\pi\cdot p/|p|$.  \emph{Setting
$\gamma\equiv b/a$ and using the parameters $a$ and $\gamma$ instead
of $a$ and $b$ codes the generic leaves so they fit smoothly into the
nongeneric ones, thus allowing the possibility of an effective
perturbation approach.}  The map $p=aw$, $\pi=\dot w+\gamma w$ for
$a=0$ is many-to-one and so the three dimensional set of nongeneric
leaves $p=0$ is ``blown-up''  by this map to the five dimensional set of
$TS^2\times\mbb R=\sset{\bigl((w,\dot w),\gamma\bigr)}$. Thus one is
led to define the \emph{blown-up space} of $\mtl{se}(3)^*$ as
\begin{equation*}
\hat P\equiv TS^2\times\mbb R_{\ge 0}\times\mbb R
\equiv\set{(w,\dot w,a,\gamma)}{|w|=1,w\cdot w=0, a\ge0}
\end{equation*}
with \emph{blow-down map}
\begin{equation*}
p=aw,\quad\pi=\dot w+\gamma w
\end{equation*}
and corresponding \emph{blow-up map}, defined on the generic ($p\ne0$)
leaves only,
\begin{equation*}
w=\frac p{|p|},\quad\dot w=\pi-\frac{\pi\cdot p}{|p|^2}p,\quad a=|p|,
\quad \gamma=\frac ba.
\end{equation*}
The blow-up map is a diffeomorphism from the (open) set of generic
leaves to the (open) set $a>0$ in the blown-up space (the
\emph{generic sector}), such that each generic leaf is sent to the
constant manifold $TS^2$.  The evolution of the generic leaves is
transformed to a evolution on $TS^2$ parameterized by the Casimir
values $a$ and $\gamma$. The blow-down map takes the set $a=0$ in the
blown-up space (the \emph{nongeneric sector}) to the set of nongeneric
leaves, and is many-to-one on that sector.  Increasing the parameter
$a$ from zero corresponds to leaving the nongeneric leaves and moving
to the generic ones, while $\gamma$ parameterizes the possible avenues
of departure.

The utility of the blow-up to support perturbation arguments depends
on whether or not the dynamics of the generic sector can be continued
smoothly to the nongeneric sector. On the generic sector the vector
field that generates the dynamics is
\begin{equation}\elb{8}\begin{split}
\frac{dw}{dt}&=\frac1a\frac{dp}{dt}=\frac1ap\times\mbf I^{-1}\pi=
w\times\mbf I^{-1}(\dot w+\gamma w),\\
\frac{d\dot w}{dt}&=\frac{d\pi}{dt}-\gamma\frac{dw}{dt}\\
&=\pi\times\mbf I^{-1}\pi+p\times\mbf M^{-1}p-\gamma\frac{dw}{dt}\\
&=(\dot w+\gamma w)\times\mbf I^{-1}(\dot w+\gamma w)+
aw\times\mbf M^{-1}aw-\gamma w\times\mbf I^{-1}(\dot w+\gamma w)\\
&=\dot w\times\mbf I^{-1}(\dot w+\gamma w)+a^2w\times\mbf M^{-1}w.
\end{split}\end{equation}
Obviously this is smooth in $a$ for all $a\ge 0$, as required.
Dynamics on $a=0$ that is robust enough to continue through
perturbation to small positive $a$ will have implications for the
original system. By continuity in $a$, the blown-up vector field is a
lift by the smooth blow-down map of the vector field for the original
system, even through the nongeneric sector. Thus the flow on the
nongeneric sector corresponds through the blow-down map to the flow of
the original system on the union of the nongeneric leaves.

Actually, the blow-up has a very transparent reformulation, since
$\hat P$ is diffeomorphic to $S^2\times\mbb R^3\times\mbb R_{\ge
0}=\sset{(w,\pi,a)}$ by the map $\pi=\gamma w+\dot w$. Through this
diffeomorphism the blow-down map is simply $p=aw$, which is to say
that $w$ by itself is enough to desingularize the foliation by
symplectic leaves, but not enough to normalize the leaves. The
blow up map is a proper map since the map $(a,w)\mapsto aw$ is proper.

Some exploration of the nongeneric sector of the blown-up space may be
helpful for visualization purposes. For fixed $\pi_0$ the equation
$\gamma w+\dot w=\pi_0$ has solution $w\in S^2$ and $\gamma=\pi_0\cdot
w$. Consequently the blow-up of the point $p=0$, $\pi=\pi_0$ is in all
cases a two sphere. This sphere intersects fixed $\gamma$ such that
$|\gamma|<|\pi_0|$ in a circle, $\gamma=\pm|\pi_0|$ in a point, and
$|\gamma|>|\pi_0|$ not at all. Thus departure from the point $p=0$,
$\pi=\pi_0$ along $\gamma>|\pi_0|$ is impossible.  The nongeneric
symplectic symplectic leaf $|\pi|=r>0$ blows up to $\gamma^2+|\dot
w|^2=r^2$ which is diffeomorphic to $S^2\times S^2$, and which for
fixed $|\gamma|<r$ is a circle bundle and for $|\gamma|=r$ is a
sphere.

The Hamiltonian pulls back through the smooth blow-down map to
\begin{equation*}\begin{split}
&\hat H\equiv\frac12\pi^t\mbf I^{-1}\pi+\frac12p^t\mbf M^{-1}p
=\hat H^0+a^2\hat H^1,\\
&\hat H^0\equiv=\frac12(\dot w+\gamma w)^t\mbf I^{-1}(\dot w +\gamma w),\quad
\hat H^1\equiv\frac{1}2w^t\mbf M^{-1}w,
\end{split}\end{equation*}
$\hat H$ is written this way in anticipation of perturbation arguments
from $a=0$ to small nonzero $a$.  The symplectic form $\hat\omega$ on
the nongeneric sector can be calculated from the formula for the
coadjoint orbit symplectic forms of $\mtl{SE}(3)$
in~\ct{MarsdenJERatiuTS-1994.1}, with the result that
\begin{equation*}\begin{split}
&\hat\omega(w,\dot w)\bigl((\delta w_1,\delta\dot w_1),
 (\delta w_2,\delta\dot w_2)\bigr)\\
&\qquad\mbox{}=-w\cdot(\delta w_1\times\delta\dot w_2
 -\delta w_2\times\delta\dot w_1)
 -\gamma w\cdot(\delta w_1\times\delta w_2).
\end{split}\end{equation*}
By continuity, the relation $i_{X_{\hat H}}\hat\omega=d\hat H$
persists from $a=0$ to $a>0$, so the vector field~\erf{8} is
Hamiltonian at $a=0$ with symplectic form $\hat\omega$ and Hamiltonian
$\hat H^0$. Thus the evolution on the nongeneric sector is
Hamiltonian, in a way that smoothly continues the Hamiltonian
structure of the generic sector.

The dynamics on the invariant submanifold $p=0$ in the original space
$P$ admits the subcasimir $|\pi|$. This conserved quantity
(\emph{conserved on $p=0$ only}) pulls back to a conserved quantity
$|w+\gamma\dot w|$ for the \emph{nongeneric sector} of the blow-up
space $\hat P$. Since $w\cdot\dot w=0$ and $|w|=1$, this gives the
conserved quantity $|\dot w|^2$ and hence the conserved quantity
$f(|\dot w|)$, where $f$ is any function. The Hamiltonian
vector field of $f(|\dot w|^2)$ is
\begin{equation*}
\frac{dw}{dt}=-\frac{f^\prime(|\dot w|)}{|\dot w|}\dot w\times w,\quad
\frac{d\dot w}{dt}=-\gamma\frac{f^\prime(|\dot w|)}{|\dot w|}w\times\dot w.
\end{equation*}
Note that $\tilde m\equiv\dot w+\gamma w$ is conserved by these
equations, so that
\begin{equation*}
\frac{dw}{dt}=-\frac{f^\prime(|\dot w|)}{|\dot w|}\tilde m\times w,\quad
\frac{d\dot w}{dt}=-\gamma\frac{f^\prime(|\dot w|)}{|\dot w|}
 \tilde m\times\dot w,
\end{equation*}
the solution of which is rotations about $\tilde m$. To normalize the
period at $2\pi$ and the righthand sense about $\tilde m$, choose
\begin{equation*}
\frac{f^\prime(|\dot w|)}{|\dot w|}|\tilde m|=
 \frac{f^\prime(|\dot w|)}{|\dot w|}\sqrt{\gamma^2+|\dot w|^2}=-1,
\end{equation*}
which gives $f(|\dot w|)=-\sqrt{\gamma^2+|\dot w|^2}$. \emph{Thus the
nongeneric sector has an additional $\mtl{SO}(2)$ symmetry}, which acts
by
\begin{equation*}
\theta\cdot(w,\dot w)\equiv\bigl(\onm{exp}(m^\wedge\theta)w,
 \onm{exp}(m^\wedge\theta)\dot w\bigr),
 \quad m\equiv\frac{\dot w+\gamma w}{\sqrt{\gamma^2+|\dot w|^2}},
\end{equation*}
and has momentum
\begin{equation*}
\hat J\equiv-\sqrt{\gamma^2+|\dot w|^2}.
\end{equation*}
This extra $\mtl{SO}(2)$ symmetry arises from a subcasimir of the
original system. The action and the corresponding momentum are defined
on the nongeneric sector where $\gamma$ and $\dot w$ are not both
zero. The set where $a=\gamma=0$ and $\dot w=0$ exactly corresponds
through the blow-down/up to the set where $p=0$ and $\pi=0$, so that
the action and its momentum are defined only on an open subset does
not affect the analysis of the relative equilibria $p_e^{\alpha_e}$.

Here are some aspects of the $\mtl{SO}(2)$ action and its relation
to the blow-down/up map.
\begin{enumerate} 
\item
The action is free except on the set $\dot w=0$, which is a 2-sphere
of fixed points. This two sphere is also the level $-|\gamma|$ of the
momentum $\hat J$, is a symplectic submanifold of $\hat P$, and as
such is equal to its own singular reduction.
\item
The orbit relation of the action together with the parameter $\gamma$
exactly absorb the additional phase space from blowing up the
nongeneric leaves. Indeed, for fixed $\gamma$ the blow-down map is a
quotient map for the action, and the orbit space is therefore smooth,
irrespective of the fact that the action is not free.
\item 
The blow-down map restricts to a quotient map for the (singular or
nonsingular) symplectic reduced space associated to the $\hat
J=\hat\mu$ level set. As such this reduced space is symplectomorphic
to the nongeneric leaf $p=0$, $|\pi|=-\hat\mu$.
\end{enumerate}
Only the verification of the third item in the nonsingular
($\hat\mu<-|\gamma|$) case is troublesome. For that, it is easily
verified that the map $\pi=\gamma w+\dot w$ is a quotient map for the
$\mtl{SO}(2)$ action on $\hat J^{-1}(\hat\mu)$ which has image the
sphere $TS_{-\hat\mu}^2=\set{\pi}{|\pi|=-\hat\mu}$. To pull down the
symplectic form~$\hat\omega$ by that map, first let $(\pi,\delta\pi_i)\in
TS_{-\hat\mu}^2$, $i=1,2$, and seek $(w,\dot w,\delta w_i,\delta\dot
w_i)$ such that
\begin{gather*}
|w|=1,\quad w\cdot\dot w=0,\quad w\cdot\delta w_i=0,\quad 
 \delta w_i\cdot\dot w+w\cdot\delta\dot w_i=0,\\
-\sqrt{\gamma^2+|\dot w|^2}=\hat\mu,\quad \delta\dot w_i\cdot\dot w=0,\quad
\pi=\gamma w+\dot w,\quad \delta\pi_i=\gamma\delta w_i+\delta\dot w_i.
\end{gather*}
To solve these equations choose and $w$ such that
$w\cdot\dot\pi=\gamma$ and set $\dot w=\pi-\gamma w$. Expanding $\delta w_i$
and $\delta\dot w_i$ in the basis $w$, $\dot w$, $w\times\dot w$ gives
\begin{equation*}
\delta w_i=-\frac{w\cdot\delta\pi_i}{\hat\mu^2-\gamma^2}(\pi-\gamma w),\quad
 \delta\dot w_i=(w\cdot\delta\pi_i)w
 +\frac{(w\times\pi)\cdot\delta\pi_i}{\hat\mu^2-\gamma^2}w\times\pi.
\end{equation*}
Substitution into~$\hat\omega$ then gives
\begin{equation*}\begin{split}
&\hat\omega(w,\dot w)\bigl((\delta w_1,\delta\dot w_1),
 (\delta w_2,\delta\dot w_2)\bigr)\\
&\qquad\mbox{}=\frac1{\hat\mu^2-\gamma^2}
 \Bigl(\bigl(\delta\pi_1\cdot w\bigr)\bigl((\pi\times\delta\pi_2)\cdot w\bigr)
 -\bigl(\delta\pi_2\cdot w\bigr)\bigl((\pi\times\delta\pi_1)\cdot w\bigr)\Bigr)
\\
&\qquad\mbox{}=\frac1{\hat\mu^2-\gamma^2}\bigl(w\times(w\times\pi)\bigr)\cdot
 (\delta\pi_1\times\delta\pi_2)\\
&\qquad\mbox{}=\frac1{\hat\mu^2-\gamma^2}\bigl((w\cdot\pi)w-\pi\bigr)\cdot
 (\delta\pi_1\times\delta\pi_2)\\
&\qquad\mbox{}=\frac1{\hat\mu^2-\gamma^2}\left((w\cdot\pi)
\frac{w\cdot\pi}{|\pi|^2}\pi
-\pi\right)\cdot
 (\delta\pi_1\times\delta\pi_2)\\
&\qquad\mbox{}=\frac1{\hat\mu^2-\gamma^2}\left(\frac{\gamma^2}{\hat\mu^2}-1
\right)\pi\cdot(\delta\pi_1\times\delta\pi_2)\\
&\qquad\mbox{}=-\frac{1}{|\pi|^2}\,\pi\cdot(\delta\pi_1\times\delta\pi_2).
\end{split}\end{equation*}
This last expression is the symplectic form on the nongeneric leaf
$|\pi|=-\hat\mu$, as required. 

In short, the $\mtl{SO}(2)$ symmetry arises and exactly absorbs the
additional dimensions resulting from the blow-up construction.  The
symplectic reductions of the nongeneric sector by this symmetry
exactly coincide with the original system restricted to the nongeneric
symplectic leaves of the phase space $P$.

The pull-back of the relative equilibria~$p_e^{\alpha_e}$ by the
blow-down map is the set of $(w,\dot w,a,\gamma)$ such that
\begin{equation*}
p=aw=0,\quad \pi=\alpha_e\mbf k=\dot w+\gamma w.
\end{equation*}
Since $|w|=1$ the first equation is equivalent to $a=0$ (the relative
equilibria are of course in the nongeneric sector), and dotting the
second with $w$ shows it is equivalent to $\gamma=\alpha_e\mbf k\cdot
w$ and $\dot w=\alpha_e\mbf k-\gamma w$. Since $\gamma^2+|\dot
w|^2={\alpha_e}^2$ and $\alpha_e\ne0$, all these solutions are in are
within the open set where the $\mtl{SO}(2)$ action is defined (i.e.
where $\gamma$ and $\dot w$ are not both zero). Thus the relative
equilibria~$p_e^{\alpha_e}$ blow up to
\begin{equation}\elb{13}
\hat p_e^{\alpha_e}:\quad|w|=1,\quad \dot w=\alpha_e\mbf k
 -\alpha_e(\mathbf k\cdot w)w,\quad a=0,\quad
 \gamma=\alpha_e\mbf k\cdot w.
\end{equation}
Since $w$ is unconstrained in~\erf{13}, except for the first equation,
$\hat p_e^{\alpha_e}$ is diffeomorphic to $S^2$. Since $m=\mbf k$ on
$\hat p_e^{\alpha_e}$, the $\mtl{SO}(2)$ action on $\hat
p_e^{\alpha_e}$ is by rotation of the pair $(w,\dot w)$ about $\mbf
k$. By substitution of $\hat p_e^{\alpha_e}$ into~\erf{8}, each point
of $\hat p_e^{\alpha_e}$ is a relative equilibrium for the
$\mtl{SO}(2)$ symmetry, except for the the two points $w=\pm\mbf k,
\dot w=0,a=0,\gamma=\pm\alpha_e$, which are equilibria that reside at
singular points of the action. Each relative equilibrium in $\hat
p_e^{\alpha_e}$ has the same generator, namely
\begin{equation}\elb{19}
\hat\xi_e^{\alpha_e}\equiv-\frac{\alpha_e}{I_3}
\end{equation}
and by substitution into $\hat J$, the same
momentum, namely 
\begin{equation*} 
\hat\mu_e^{\alpha_e}\equiv-\alpha_e.
\end{equation*}

Fixing $\gamma$, which means fixing a parameter,
$\alpha_e=\gamma/\cos\phi$, where $\phi$ is the angle between $\mbf k$
and $w$. Thus for fixed $\gamma$ there are 2-submanifolds of relative
equilibria, as $\alpha_e$ is varied, as expected for an $\mtl{SO}(2)$
symmetric Hamiltonian system.  Along those submanifolds there is the
\emph{momentum-generator} relation
\begin{equation}\elb{22}
\hat\xi_e^{\alpha_e}=\frac1{I_3}\hat\mu_e^{\alpha_e},
\end{equation}
which will give a crucial component in the KAM twist condition to
follow.

\section{Normal forms in the blown-up system}
Proposition~\thrf{StabProp} below relates the stability of the
relative equilibrium~$p_e^{\alpha_e}$ to the stability of its blow-up
$\hat p_e^{\alpha_e}$.

\begin{proposition}\lb{StabProp}
Suppose that, for some fixed $\alpha_e$, $\hat p_e^{\alpha_e}$ is
stable for the flow $\hat F_t$ on $\hat P$, in the sense that for all
neighborhoods $\hat U$ of $\hat p_e^{\alpha_e}$ there is a
neighborhood $\hat V$ of $\hat p_e^{\alpha_e}$ such that $\hat
F_t(\hat p)\in\hat U$ for all $\hat p\in\hat V$. Then $p_e^{\alpha_e}$
is a stable relative equilibrium.
\end{proposition}

\begin{proof}
Suppose $U$ is a neighborhood of $p_e^{\alpha_e}$. $U$ pulls back by
the blow-down map to an open neighborhood $\hat U$ of $\hat
p_e^{\alpha_e}$. Let $\hat V$ be a neighborhood as in the statement of
the proposition. Then it suffices to show that $\hat V$ pushes forward
by the blow-down map to a neighborhood of $p_e^{\alpha_e}$. But this
follows since the blow-down map is proper.
\end{proof}

In particular, if all of the relative equilibria and both equilibria
in $\hat p_e^{\alpha_e}$ are stable under perturbation
both within the phase space $TS^2$ and in the parameters $a$ and
$\gamma$, then the original relative equilibrium~$p_e^{\alpha_e}$ is
stable. When $a$ is perturbed away from $0$ this is a $\mtl{SO}(2)$
symmetry breaking perturbation. As $\onm{dim}TS^2=4$, the blown-up
system is integrable when $a=0$ and hence the stability issue is one
of the stability of periodic orbits of a near integrable Hamiltonian
system.

Assume, without loss of generality, that $\alpha_e>0$.  Since the
$\mtl{SO}(2)$ symmetry on $\hat p_e^{\alpha_e}$ is by  rotation about
$\mbf k$, it suffices to consider the stability of orbits in $\hat
p_e^{\alpha_e}$ emanating from points $\hat p_e^{\alpha_e,\theta}$
obtained by substituting $w=\sin\theta\mbf i+\cos\theta\mbf k$
into~\erf{13}, for $\theta\in[0,\pi]$.

\subsection{Normal form for the relative equilibria}\lb{200}
Consider first the \emph{relative equilibria} in $\hat
p_e^{\alpha_e,\theta}$; i.e., exclude the \emph{equilibria}
corresponding to $\theta=0$ and $\theta=\pi$.  The argument proceeds
by adapting and incrementally refining, to the order required for the
stability analysis, the normal form near relative equilibria developed
in~\ct{PatrickGW-1995.1}.

Below $O(x;y)^k$, $x\in\mbb R^n$, $y\in\mbb R^m$ will denote the set
of smooth $y$ dependent functions such that $O(x;y)/|x|^k$ is bounded
near $0$. The product $O(x;y)^kO(x^\prime;y^\prime)^{k^\prime}$ denotes
the set of finite sums of products of elements of $O(x;y)^k$ and
$O(x^\prime;y^\prime)^{k^\prime}$.

\subsubsection{Initial normal form}\lb{201}
This is constructed from the linearization of the relative
equilibrium, which means the linearization at $\hat p_e^{\alpha_e,\theta}$
of Hamiltonian vector field $X_{\hat H_{\xi^{\alpha_e}_e}}$ where
\begin{equation*}
\hat H_{\xi^{\alpha_e}_e}^0\equiv\hat H^0-\hat\xi_e\hat J.
\end{equation*}
The characteristic polynomial of the linearization is $x\mapsto
x^2(x^2+{\omega_e}^2)$, where
\begin{equation}\elb{23}
\omega_e\equiv\pm{\alpha_e}\sqrt{\left(\frac1{I_3}-\frac1{I_1}\right)
 \left(\frac1{I_3}-\frac1{I_2}\right)},
\end{equation}
For later convenience define $\omega_e$ positive if $I_3>I_1$ and
$I_3>I_2$ and negative if $I_3<I_1$ and $I_3<I_2$.  The linearization
has a $0$ and $\pm i\omega_e$~generalized eigenspaces, both of
dimension~$2$.  Introducing the parameter
\begin{equation*}
D\equiv\frac{I_2(I_3-I_1)}{I_1(I_3-I_2)},
\end{equation*}
the vectors
\begin{equation*}\begin{split}
&v_1\equiv\frac {D^{\frac14}}{\sqrt{\alpha_e}}\left[\begin{array}{cccccc}0&\cos\theta&0&0&\alpha_e\sin^2\theta&0
 \end{array}\right]\\
&v_2\equiv\frac {D^{-\frac14}}{\sqrt{\alpha_e}}\left[
\begin{array}{cccccc}\cos\theta&0&-\sin\theta&\alpha_e\sin^2\theta&0&
  \alpha_e\sin\theta\cos\theta\end{array}\right]\\
&v_3\equiv\sin\theta\left[
\begin{array}{cccccc}0&1&0&0&-\alpha_e\cos\theta&0
 \end{array}\right]\\
&v_4\equiv\frac1{\alpha_e\sin\theta}\left[
 \begin{array}{cccccc}-\cos^2\theta&0&\cos\theta\sin\theta
 &\alpha_e\cos^3\theta&0&-\alpha_e\sin\theta(1+\cos^2\theta)\end{array}\right]
\end{split}\end{equation*}
form a basis of $T_{\hat p_e^{\alpha_e,\theta}}S^2$ which satisfies
the following:
\begin{enumerate}
\item the basis is symplectically canonical, so that the symplectic form
with respect to it is
\begin{equation*}
\omega(\hat p_e^{\alpha_e,\theta})=\left[
\begin{array}{cccc}0&1&0&0\\-1&0&0&0\\0&0&0&1\\0&0&-1&0\end{array}\right];
\end{equation*}
\item with respect to the basis the derivative of the momentum is
\begin{equation*}
d\hat J\bigl(\hat p_e^{\alpha_e,\theta}\bigr)=
 \left[\begin{array}{cccc}0&0&0&1\end{array}\right];
\end{equation*}
\item the third basis vector~$v_3$ is the infinitesimal generator
action corresponding to~$1\in\mtl{so}(2)$;
\item the first two vectors $v_1$, $v_2$
span the~$\omega_e$ generalized eigenspace and the last two $v_3$, $v_4$ span
the $0$~generalized eigenspace of the linearization;
\item the linearization of the relative equilibrium is
\begin{equation}\elb{27}
dX_{\hat H_{\xi^{\alpha_e}_e}^0}(\hat p_e^{\alpha_e,\theta})
  =\left[\begin{array}{cccc}
  0&\omega_e&0&0\\-\omega_e&0&0&0\\0&0&0&\kappa_e\\0&0&0&0\end{array}\right],
\end{equation}
where
\begin{equation}\elb{25}
\kappa_e\equiv 1/I_3.
\end{equation}
\end{enumerate}
Consequently, the basis effects a Witt-Moncrief decomposition
\begin{equation*}
T_{\hat p_e^{\alpha_e,\theta}}\hat P=N_1\oplus\mtl{so}(2)
\oplus\mtl{so}(2)^*
\equiv\onm{span}(v_1,v_2)\oplus\mbb Rv_3\oplus\mbb Rv_4.
\end{equation*}
Here $N_1$, the \emph{symplectic normal}, may be
identified with the tangent space to the symplectic reduced space at
$\hat p_e^{\alpha_e,\theta}$ for the $\mtl{SO}(2)$ action. The
appearance of the nilpotent part of the linearization is the
foundational element of~\ct{PatrickGW-1995.1}. The value of
$\kappa_e$ coincides with the derivative
$d\hat\xi_e^{\alpha_e}/d\hat\mu_e^{\alpha_e}$
from~\erf{22}, as predicted by the general theory. 

The initial normal form can now transcribed from the data above, and is
\begin{equation}\begin{split}\elb{21}
&\hat H=\frac{\omega_e}2(q^2+p^2)+\hat\xi_e^{\alpha_e}\nu
 +\frac12\kappa_e\nu^2+R(q,p,\nu)
 +\frac{a^2}{2}\hat H^1(q,p,\varphi,\nu),\\
&R=O(q,p,\nu)^3
\end{split}
\end{equation}
on the product of $\mbb R^2\times
T^*\mtl{SO}(2)=\sset{(q,p),(\varphi,\nu)}$ with the product symplectic
form $dq\wedge dp+d\varphi\wedge d\nu$, with $\mtl{SO}(2)$ acting by
lifts of its left action on itself, and with the momentum map
$\nu-\alpha_e$.  The transcription is that there is an
$\mtl{SO}(2)$~intertwining symplectic diffeomorphism from a
neighborhood of the $\mtl{SO}(2)$ orbit of $\hat p_e^{\alpha_e}$
to a neighborhood of $0$ times the zero section of $T^*\mtl{SO}(2)$
which
\begin{enumerate}
\item sends the relative equilibrium $\hat p_e^{\alpha_e,\theta}$ to
$p=q=\nu=\varphi=0$.
\item has derivative at $\hat p_e^{\alpha_e}$ the identity map with
respect to the basis~$v_i$ and the standard basis of $\mbb R^2\times
T^*\mtl{SO}(2)$;
\item intertwines the momentum maps $\hat J$ and $\nu-\alpha_e$;
\end{enumerate}
Thus the transription is \emph{strucure preserving} in that it is
symplectic and it preserves the $\mtl{SO}(2)$ symmetry and momentum,
so the blown-up system near the group orbit of the relative
equilibrium $\hat p_e^{\alpha_e,\theta}$ can be replaced by the
entirely equivalent system~\erf{21} near $q=p=\nu=0$.

\subsubsection{Elimination of $qO(\nu)^2$, $pO(\nu)^2$, and $(q^2-p^2)\nu$}\lb{202}
The remainder term of~\erf{21} can be expanded as
\begin{equation*}\begin{split}
R=&c_1(\nu)q+c_2(\nu)p+c_3\nu(q^2-p^2)+c_4\nu(q^2+p^2)\\
&\qquad\qquad+O(q,p)^3+O(\nu)^3+O(q,p)^2O(q,p,\nu)^2,
\end{split}\end{equation*}
where $c_1(\nu)=O(\nu)^2$,$c_2(\nu)=O(\nu)^2$, and $c_3$, $c_4$ are
constants.  The transformation
\begin{equation*}
\tilde q=q+\frac{c_1}{\omega_e},\quad
  \tilde\varphi=\varphi+\frac{p\nu}{\omega_e}\,\frac{dc_1}{d\nu},
\end{equation*}
suggested by completing the square in $\frac12\omega_eq^2+c_1q\nu^2$,
is structure preserving and changes the Hamiltonian to the same form
but without terms of the form $qO(\nu)^2$. Similarly one eliminates
$pO(\nu)^2$.  The transformation
\begin{equation}\elb{64}
\tilde q=\frac q{f(\nu)},\quad\tilde p=f(\nu)p,
  \quad f(\nu)=\left(
  \frac{1-\frac{2c_3}{\omega_e}\nu}{1+\frac{2c_3}{\omega_e}\nu}
  \right)^{\frac14},
\end{equation}
takes the fragment $\frac{\omega_e}2(q^2+p^2)+c_3(q^2-p^2)\nu$ to
\begin{equation*}
\left(\frac{\omega_e}2+c_3\right)q^2+\left(\frac{\omega_e}2-c_3\right)p^2
=\frac{\omega_e}2(\tilde q^2+\tilde p^2)+\tilde q^2O(\nu)^2+\tilde p^2O(\nu)^2,
\end{equation*}
while the symplectic form becomes
\begin{equation*}\begin{split}
dq\wedge dp+d\varphi\wedge d\nu=&d\tilde q\wedge d p+d\varphi\wedge d\nu+
\frac{f^\prime}{f}(qdp+pdq)\wedge d\nu\\
=&d\tilde q\wedge d p+d\left(\varphi+\frac{f^\prime}{f}qp\right)\wedge d\nu.
\end{split}\end{equation*}
Adjoining $\tilde\varphi=\varphi+(f^\prime/f)qp$ to~\erf{64} gives a
structure preserving symplectic transformation that
eliminates the term $c_3(q^2-p^2)$.  Thus, without loss of generality,
\begin{equation}\elb{67}
R=c_4(q^2+p^2)+O(q,p)^3+O(\nu)^3+O(q,p)^2O(q,p,\nu)^2.
\end{equation}

\subsubsection{Normal form for the rigid body}\lb{203}
We will require the first two terms of the normal form corresponding
to the equilibrium $\pi=\alpha_e\mbf k$ of the blown-up system reduced
by its $\mtl{SO}(2)$ symmetry, i.e. the symplectic reduced spaces of
the rigid body $\frac12\pi\mbf I^{-1}\pi$.  The map
\begin{equation*}
\pi=\Bigl(\bigl(\alpha_e-{\textstyle\frac14}(Q^2+P^2)\bigr)^{\frac12}P,
\bigl(\alpha_e-{\textstyle\frac14}(Q^2+P^2)\bigr)^{\frac12}Q,\alpha_e
{\textstyle-\frac12}(Q^2+P^2)\Bigr)
\end{equation*}
is a symplectic chart on the reduced space $|\pi|=\alpha_e$ and in
these coordinates, the Hamiltonian becomes, up to a constant,
\begin{equation*}
\frac12\pi^t\mbf I^{-1}\pi=\frac12\bigl(\alpha_e
  -{\textstyle\frac14}(Q^2+P^2)\bigr)
 \left(\left(\frac1{I_1}-\frac{1}{I_3}\right)P^2+
 \left(\frac1{I_2}-\frac{1}{I_3}\right)Q^2\right).
\end{equation*}
Action-angle variables for the linearized flow are
\begin{equation*}
Q=\sqrt{2I}D^{\frac14}\sin\psi,\qquad P=\sqrt{2I}D^{-\frac14}\cos\psi, 
\end{equation*}
and the Hamiltonian is then
\begin{equation*}
\frac12\pi^t\mbf I^{-1}\pi=\omega_eI-\frac{\omega_e}{2\alpha_e}(D^{\frac12}\sin^2\psi
 +D^{-\frac12}\cos^2\psi)I^2.
\end{equation*}
By averaging over $\psi$,
\begin{equation}\elb{14}
\frac12\pi^t\mbf{I}\pi=\omega_eI+\frac12\upsilon_eI^2+O(Q,P)^5
\end{equation}
where
\begin{equation}\elb{24}
\upsilon_e=-\frac{\omega_e}{2\alpha_e}(D^{\frac12}+D^{-\frac12})
=\frac12\left(\frac2{I_3}-\frac1{I_1}-\frac1{I_2}\right).
\end{equation}

\subsubsection{Matching and normalizing the reduced spaces at 
$\hat p_e^{\alpha_e,\theta}$}\lb{204}
For $a=0$ the symplectic reduced space through $\hat
p_e^{\alpha_e,\theta}$ of the blown-up system is the $|\pi|=\alpha_e$
symplectic reduced space of the rigid body $\frac12\pi^t\mbf I\pi$.
For $a=0$ the symplectic reduced space of the normal form~\erf{21}
through $q=p=\varphi=\nu=0$ is $\mbb R^2$ with symplectic form
$dq\wedge dp$ and Hamiltonian $\hat H|_{\nu=0}$. Since the
intertwining map between the blown-up system and the normal form is
structure preserving, it descends to symplectomorphisms of reductions
of the these two systems. Consequently, by symplectomorphism on
$(q,p)$ only, the normal form Hamiltonian~\erf{21} at $q=p=\nu=0$ can
be equated to the rigid body normal form~\erf{14}, after which the
normal is correct to fourth order in pure $q$ and $p$ and
\begin{equation*}
R=O(q,p)^5+\nu O(q,p)^2O(q,p,\nu).
\end{equation*}

\subsubsection{Refinement by matching the normal forms and 
generators along the relative equilibria near $\hat
p_e^{\alpha_e,\theta}$}\lb{205} 
Advantage may be obtained by comparing reduced normal along the
relative equilibria $\hat p_e^{\alpha_e+z,\theta}$ as $z$
varies. These relative equilibria occur (for both systems) at momentum
$-(\alpha_e+z)$.  For the rigid body the only $z$ dependent adjustment
is in the $\alpha_e$ dependence of the linearized frequency, which
becomes $\omega_e(\alpha_e+z)/\alpha_e$, so the normal form is
\begin{equation}\elb{72}
\left(\omega_e+\frac{\omega_e}{\alpha_e}z\right)I+\frac12\upsilon_eI^2
  +O(q,p;z)^5.
\end{equation}
For~\erf{21} it is the normal form of the reduction at $\nu=-z$, so it
is the normal form of the Hamiltonian
\begin{equation*}
\omega_eI+\frac12\upsilon_eI^2-2c_4Iz+O(q,p)^2O(q,p,z)^2
\end{equation*}
which is
\begin{equation}\elb{73}
\bigl(\omega_e-2c_4z+O(z)^2\bigr)I+O(q,p;z)^3.
\end{equation}
Comparison of~\erf{72} and~\erf{73} at first order in $I$ gives a
crucial fact: 
\begin{equation*}
c_4=-\frac{\omega_e}{2\alpha_e}.
\end{equation*}

Also, the $SO(2)$ generator of the blown-up system at $\hat
p_e^{\alpha_e+z,\theta}$, which is $-(\alpha_e+z)/I_3$, and the
$\mtl{SO}(2)$ generator of system~\erf{21} at the relative equilibrium
$q=p=0$ are the same. Equating these gives
\begin{equation*}
\hat\xi_e^{\alpha_e}+\kappa_e\nu+\left.\frac{\partial R}{\partial\nu}\right|_
{\!\!\!\renewcommand{\arraystretch}{.2}\begin{array}[c]{l}\scriptstyle q=p=0\\\scriptstyle\nu=-z\end{array}}\!\!\!=-\frac1{I_3}(\alpha_e+z)=\hat\xi_e^{\alpha_e}-\kappa_ez,
\end{equation*}
which means that $R$ has no pure $\nu$ terms. Particularly, the
$O(\nu)^3$ term in~\erf{67} is zero.

\subsubsection{Symmetry breaking term}\lb{206}
The transcription to the initial normal form is known to first order since it
has derivative the identity map along the $\mtl{SO}(2)$ orbit of $\hat
p_e^{\alpha_e,\theta}$.  Consequently, $\hat H^1$ can be calculated to
first order by substitution of
\begin{equation*}
w=\exp(\varphi\mbf k)P_w\bigl(\hat p_e^{\alpha_e,\theta}+(qv_1+pv_2
+\nu v_4)\bigr)
\end{equation*}
into $w^t\mbf M^{-1}w$, where $P_w(w,\dot w)=w$.

\subsubsection{Altogether}\lb{207}
Putting all the foregoing together, the normal form is
\begin{equation}\begin{split}\elb{69}
\hat H=&\omega_eI+\frac12\upsilon_eI^2+\hat\xi_e^{\alpha_e}\nu
+\frac12\kappa_e\nu^2-\frac{\omega_e}{\alpha_e}I\nu\\
&\quad\mbox{}+O(q,p)^5+\nu O(q,p)^2O(q,p,\nu)\\
&\quad\mbox{}+a^2\hat H^{1,0}(q,p,\varphi,\nu)+a^2\hat H^{1,1}(q,p,\varphi,\nu)
+a^2O(q,p,\nu;\varphi)^2,
\end{split}\end{equation}
where
\begin{equation}
\hat H^{1,0}\equiv\frac{M_2-M_1}{2M_1M_2}\sin^2\theta\cos^2\varphi
\end{equation}
and
\begin{equation*}\begin{split}
\hat H^{1,1}&\equiv\frac{(M_2-M_1)\sin2\theta}{4M_1M_2\sqrt{\alpha_e}}
\left(-D^{\frac14}q\sin2\varphi+D^{-\frac14}p\cos2\varphi\right)\\
&\qquad\mbox{}-\frac{\cos^2\theta}{\alpha_e}\left(\frac1{M_1}\cos^2\varphi
  +\frac1{M_2}\sin^2\varphi-\frac1{M_3}\right)\nu\\
&\qquad\mbox{}-\frac{\sin2\theta}{4D^{\frac14}\sqrt\alpha_e}
  \left(\frac2{M_3}-\frac1{M_1}-\frac1{M_2}\right)p.
\end{split}\end{equation*}
The details of the symmetry breaking term $\hat H^1$ are not required
for the stability analysis and are displayed here for the sake of
completeness.  The functional form of $\hat H^{1,1}$ depends on the
choice of the basis $v_i$ and further normalization or analysis would
be required to extract information from it.

\subsection{Normal form for the equilibria}\lb{300}
There remains to consider the two equilibria $\hat p_e^{\alpha_e,0}$ and
$\hat p_e^{\alpha_e,\pi}$ corresponding to $w=\mbf k$ and $w=-\mbf k$,
respectively. These equilibria are fixed points of the action of
$\mtl{SO}(2)$ and the analysis requires a transparent extension of the
normal form in~\ct{PatrickGW-1995.1} to \emph{equilibria} which have
$\mtl{SO}(2)$ isotropy.

It suffices to consider $\hat p_e^{\alpha_e,0}$; the case of $\hat
p_e^{\alpha_e,\pi}$ is similar.  There is a one parameter family of
possible linearizations of the equilibrium, namely the linearizations
at $\hat p_e^{\alpha_e,0}$ of the Hamiltonian vector fields $X_{\hat
H_\lambda^0}$ where $\hat H_\lambda^0-\lambda\hat J$.  These
linearizations have characteristic polynomials
\begin{equation*}
x\mapsto\bigl(x^2+(\alpha_e+\lambda I_3)^2\bigr)\bigl(x^2+{\omega_e}^2\bigr).
\end{equation*}
Choosing $\lambda=-\alpha_e/I_3$ gives the largest possible null space
and therefore the largest number of intrinsically defined higher order
terms. The vectors
\begin{equation*}\begin{split}
&v_{1,0}\equiv\frac{D^{\frac14}}{\sqrt{\alpha_e}}
 \left[\begin{array}{cccccc}0&1&0&0&&0\end{array}\right],\quad
v_{2,0}\equiv\frac {D^{-\frac14}}{\sqrt{\alpha_e}}\left[
\begin{array}{cccccc}1&0&0&0&0&0\end{array}\right]\\
&v_{3,0}\equiv\frac 1{\sqrt{\alpha_e}}\left[
 \begin{array}{cccccc}0&1&0&0&-\alpha_e&0\end{array}\right]\quad
v_{4,0}\equiv\frac 1{\sqrt{\alpha_e}}\left[
\begin{array}{cccccc}-1&0&0&\alpha_e&0&0
 \end{array}\right],
\end{split}\end{equation*}
form a basis of $T_{\hat p_e^{\alpha_e,0}}S^2$ which is symplectically
canonical and such that the linearization
$dX_{H_\lambda^0}(\hat p_e^{\alpha_e,0})$ is the
linearization~\erf{27} with $\kappa_e$ replaced by zero.  The first
two vectors span the tangent space to the (singular) reduced space
through $p_e^{\alpha_e,0}$, which is the $|\pi|=\alpha_e$ symplectic
reduced space of the rigid body $\frac12\pi^t\mbf I^{-1}\pi$.

Higher order terms of the Taylor expansion of $X_{H_\lambda^0}$ are
intrinsically polynomials on $\onm{ker}dX_{H_\lambda^0}(\hat
p_e^{\alpha_e,0})$.  Alternately one can compute the higher order
terms of the Taylor expansion of the Hamiltonian on the null space.
Letting $(x,y)$ be the coordinates on $\onm{ker}dX_{H_\lambda^0}(\hat
p_e^{\alpha_e,0})$ indicated by the last of the two basis vectors
above, the Hamiltonian on the null space is easily computed to be
\begin{equation*}
\frac{1}{8I_3}(x^2+y^2)^2+O(x,y)^5=\frac{\kappa_e}8(x^2+y^2)^2+O(x,y)^5.
\end{equation*}
The initial normal form, obtained by the Equivariant Darboux Theorem,
is the Hamiltonian
\begin{equation*}\begin{split}
\hat H=&\frac{\omega_e}2(q^2+p^2)
  -\frac{\alpha_e}{2I_3}(x^2+y^2)+\frac{\kappa_e}{8}(x^2+y^2)^2\\
&\qquad\mbox{}+O(q,p)O(q,p,x,y)^2+O(x,y)^5+a^2\hat H^1(q,p,x,y)
\end{split}\end{equation*}
on the phase space $\mbb R^2\times\mbb R^2=\sset{(q,p),(x,y)}$, with
symplectic form $dq\wedge dp+dx\wedge dy$, with $\mtl{SO}(2)$ acting
by counterclockwise rotation on $(x,y)$, and with the momentum mapping
$-\frac12(x^2+y^2)$. The transcription is by local symplectic
diffeomorphism with analogous properties to those stated in
Section~\thrf{201}.

Manipulations similar to those in Section~\thrf{200} are required, as follows:
\begin{enumerate}
\item Linear terms in $q$ of the form $qO(x,y)^2$ can be removed as in
Section~\thrf{202}, and similarly linear terms in $p$, $x$, and $y$.
Quartic terms of the form $O(q,p)^2O(x,y)^2$ must  by $\mtl{SO}(2)$ invariance
be in $(x^2+y^2)O(q,p)$, and so can be written as sums of
$(q^2+p^2)(x^2+y^2)$ and $(q^2-p^2)(x^2+y^2)$, and the latter kind removed,
as in Section~\thrf{202}.
\item Pure $q$ and $p$ terms up to order~4 can be found by matching the
reduced system of the initial normal form to  rigid body reduced spaces.
\item By Items~1 and~2 all terms up to and including order 4 are
removed or calculated, except for the coefficient of the term
$(q^2+p^2)(x^2+y^2)$. This can be found by matching normal forms along
the equilibria $q=p=0$ (which are fixed points of the action of
$\mtl{SO}(2)$), and the resulting term is
$-\frac{\omega_e}{4\alpha_e}(q^2+p^2)(x^2+y^2)$.
\item The remainder after all of that, having no terms linear in any
variable, and being at least degree $5$, is of the form
$O(q,p,x,y)^2O(q,p,x,y)^3$, and is $\mtl{SO}(2)$ invariant.
\item The symmetry breaking term $\hat H^1$ can be calculated as in 
Section~\thrf{206} by substituting
\begin{equation*}
w=P_{w}(\hat p_e^{\alpha_e,0}+qv_{1,0}+pv_{2,0}+xv_{3,0}+yv_{4,0})
\end{equation*}
into $w^t\mbf M^{-1}w$, and keeping the leading terms, which are order~2.
\end{enumerate}
Altogether, the normal form is
\begin{equation}\begin{split}\elb{325}
\hat  H=&\omega_eI+\frac12\upsilon_eI^2
  +\xi_e^{\alpha_e}\nu+\frac12\kappa_e\nu^2-\frac{\omega_e}{\alpha_e}I\nu
 +O(q,p,x,y)^2O(q,p,x,y)^3\\
&\qquad\qquad\mbox{}+\frac{a^2}2\hat H^{1,1}(q,p,x,y)+a^2O(q,p,x,y)^3
\end{split}\end{equation}
where $I=\frac12(q^2+p^2)$ and $\nu=\frac12(x^2+y^2)$, and
\begin{equation*}
\hat H^{1,1}\equiv\frac{(M_3-M_1)(D^{-\frac14}p-y)^2}{\alpha_eM_1M_3}+
\frac{(M_3-M_2)(D^{\frac14}q+x)^2}{\alpha_eM_2M_3}.
\end{equation*}

\section{Stability}
The rescaling
\begin{equation*}
I=a^{2c}\tilde I,\quad\nu=a^{2c}\tilde\nu
\end{equation*}
is symplectic with multiplier $a^{2c}$. Substituting into~\erf{69} and 
dropping the tildes gives
\begin{equation*}\begin{split}
\frac1{a^{2c}}\hat H&=\omega_eI+\frac12\upsilon_eI^2a^{2c}
+\hat\xi_e^{\alpha_e}\nu+\frac12\kappa_e\nu^2a^{2c}-\omega_eI\nu a^{2c}
+O(a;q,p,\varphi,\nu)^{3c}\\
&\qquad\mbox{}+a^{2-2c}\hat H^{1,0}+O(a;q,p,\varphi,\nu)^{2-c}
\end{split}\end{equation*}
Matching the exponents of $a$ in first nontrivial terms of the
integrable part, i.e. $\frac12\kappa_e\nu^2a^{2c}$ and
$\frac12\upsilon_eI^2a^{2c}$, with the first term of the nonintegrable
part, gives $2c=2-2c$, or $c=\frac12$. After putting $\epsilon=\sqrt
a$, and disposing the factor $1/a^{2c}$ of $\hat H$, which merely
reparameterizes time, one has
\begin{equation}\elb{15}
\hat H=\omega_eI+\hat\xi_e^{\alpha_e}\nu+\left(\frac12\kappa_e\nu^2
  -\frac{\omega_e}{\alpha_e}I\nu+\frac12\upsilon_eI^2+\hat H^{1,0}\right)\epsilon^2+O(\epsilon)^3,
\end{equation}
where the dependence of $O(\epsilon)^3$ on all of $q$, $p$, $\varphi$,
and $\nu$ has been notationally suppressed.  For $\epsilon=0$ the
Hamiltonian~\erf{15} has a periodic orbit cylinder by varying
$\varphi$ and $\nu$ with $I=0$. The orbit $I=\nu=0$ corresponds to the
relative equilibrium $\hat p_e^{\alpha_e,\theta}$.  

For determining stability it suffices to approximate the Poincar\'e
map for the orbit corresponding to $\nu=0$ in the zero energy
level. Solving~\erf{15} for $\nu$ when $H=0$ gives
\begin{equation*}
\nu=-\frac{\omega_eI}{\hat\xi_e^{\alpha_e}}+O(\epsilon)^2
\end{equation*}
and the equations of motion for~\erf{15} are
\begin{equation}\elb{16}\begin{split}
&\frac{d\psi}{dt}=\omega_e+\left(\upsilon_eI-\frac{\omega_e}{\alpha_e}\nu\right)\epsilon^2
  +O(\epsilon)^3,
\quad\frac{dI}{dt}=O(\epsilon)^3,\\
&\frac{d\varphi}{dt}=\hat\xi_e^{\alpha_e}+\left(\kappa_e\nu-\frac{\omega_e}{\alpha_e}I\right)\epsilon^2
  +O(\epsilon)^3,
\quad\frac{d\nu}{dt}=-\frac{\partial\hat H^{1,0}}{\partial\varphi}\epsilon^2
 +O(\epsilon)^3.
\end{split}\end{equation}
On the zero energy level the equation for the evolution of $\varphi$ is
\begin{equation}\elb{17}
\frac{d\varphi}{dt}=\hat\xi_e^{\alpha_e}
  -\frac{\omega_e(\kappa_e\alpha_e+\hat\xi_e^{\alpha_e})}{\alpha_e\hat\xi_e^{\alpha_e}}
  I\epsilon^2+O(\epsilon)^3.
\end{equation}
With an initial condition $I=I_0$, the second equation of~\erf{16}
gives $I=I_0+O(\epsilon)^3$, so the return time of the Poincar\'e map
is, from~\erf{17},
\begin{equation*}
T\equiv-\frac{2\pi}{\hat\xi_e^{\alpha_e}}\left(1+\frac{(\alpha_e\kappa_e
  +\hat\xi_e^{\alpha_e})\omega_eI_0}{\alpha_e(\hat\xi_e^{\alpha_e})^2}
  \epsilon^2\right)+O(\epsilon)^3
\end{equation*}
Solving the first two equations of~\erf{16} over this period,
dropping the subscript~0 for the initial conditions, and using
the relation $\hat\xi_e^{\alpha_e}=-\kappa_e\alpha_e$ to eliminate $\alpha_e$
gives
\begin{equation}\begin{split}\elb{18}
&I^\prime=I+O(\epsilon)^3\\
&\psi^\prime=
\psi-\frac{2\pi\omega_e}{\hat\xi_e^{\alpha_e}}
-\epsilon^2\frac{2\pi}{(\hat\xi_e^{\alpha_e})^3}
  \bigl(\upsilon_e(\hat\xi_e^{\alpha_e})^2-\kappa_e{\omega_e}^2\bigr)I
  +O(\epsilon)^3
\end{split}\end{equation}
This is of the form (\ct{MeyerKRHallGR-1991.1}, Theorem~2, page~231),
namely $(I,\psi)\mapsto(I^\prime,\psi^\prime)$ by
\begin{equation*}\begin{split}
&I^\prime=I+\epsilon^{r+s}c(I,\psi,\epsilon),\\
&\psi^\prime=\psi+\omega+\epsilon^sh(I)+\epsilon^{s+r}d(I,\psi,\epsilon)
\end{split}\end{equation*}
with $r=1$, $s=2$, and
\begin{equation*}
h(I)\equiv-\frac{2\pi}{(\hat\xi_e^{\alpha_e})^3}
\bigl(\upsilon_e(\hat\xi_e^{\alpha_e})^2-\kappa_e{\omega_e}^2\bigr)I.
\end{equation*}
The twist condition $dh/dI\ne0$ is
\begin{equation}\elb{30}
\upsilon_e(\hat\xi_e^{\alpha_e})^2-\kappa_e{\omega_e}^2\ne0,
\end{equation}
which, after substituting~\erf{19}, \erf{23}, \erf{25}, and~\erf{24},
is 
\begin{equation*}\begin{split}
&-\frac12\left(\frac2{I_3}-\frac1{I_1}-\frac1{I_2}\right)+I_3\left(\frac1{I_3}-\frac1{I_1}\right)\left(\frac1{I_3}-\frac1{I_2}\right)\\
&\qquad\qquad\qquad\mbox{}=\frac{1}{I_1I_2}\left(I_3-\frac12(I_1+I_2)\right)\ne0.
\end{split}\end{equation*}
This is certainly true  if $I_3$ is not  between $I_1$ and $I_2$.

As for the equilibria $\hat p_e^{\alpha_e,0}$ of Section~\thrf{300},
Arnold's Stability Theorem (\ct{MeyerKRHallGR-1991.1}, Theorem~1,
page~235) together with the normal form~\erf{325} imply that the
equilibrium is stable when $a=0$ if
\begin{equation*}
\left.\omega_eI+\frac12\upsilon_eI^2
  +\hat\xi_e^{\alpha_e}\nu+\frac12\kappa_e\nu^2-\frac{\omega_e}{\alpha_e}I\nu
  \right|_
{\!\!\!\renewcommand{\arraystretch}{.2}\begin{array}[c]{l}\scriptstyle
 I=\hat\xi_e^{\alpha_e}\\\scriptstyle\nu=-\omega_e\end{array}}
\!\!\!\!\!\!=\frac12(\upsilon_e(\hat\xi_e^{\alpha_e})^2-\kappa_e{\omega_e}^2)\ne0,
\end{equation*}
which is the same as the twist condition~\erf{30}.  Stability follows
for sufficiently small nonzero~$a$ since $a$ contributes continuously.

Thus $p_e^{\alpha_e}$ are stable as equilibria on the Poisson reduced
space $\sset{(\pi,p)}$. From~\ct{PatrickGWRobertsRMWulffC-2001.1}, and
since the momentum at $p_e^{\alpha_e}$ has zero translational part and
rotational part parallel to $\mbf k$, this implies $\mtl{SO}(2)\times
C$ stability of the original equilibrium, where $C$ is any cone about
$\mbf k$.

\begin{figure}[t]
\setlength{\unitlength}{1in}
\centerline{{\begin{picture}(4.8,1.225)
\put(.1,.04){{\epsfysize=1.2in\epsfbox{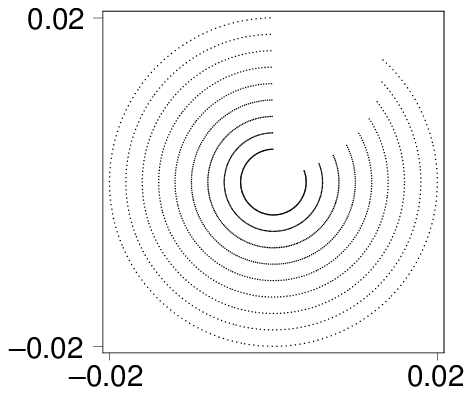}}}
\put(0,.9){\footnotesize$\sqrt[4]{D}\pi_x$}
\put(.8,0){\footnotesize$\pi_y/\sqrt[4]{D}$}
\put(1.65,.04){{\epsfysize=1.2in\epsfbox{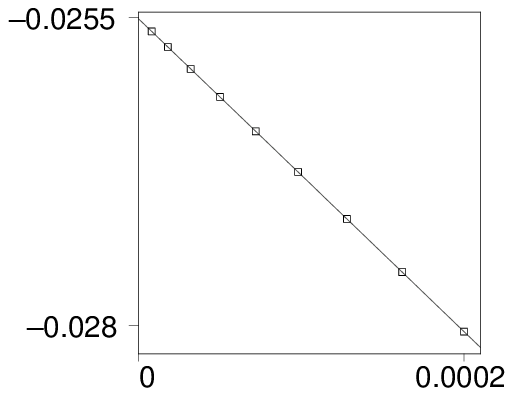}}}
\put(1.64,.9){\footnotesize$\psi^\prime-\psi$}
\put(2.5,0){\footnotesize$I$}
\put(3.425,.04){{\epsfysize=1.2in\epsfbox{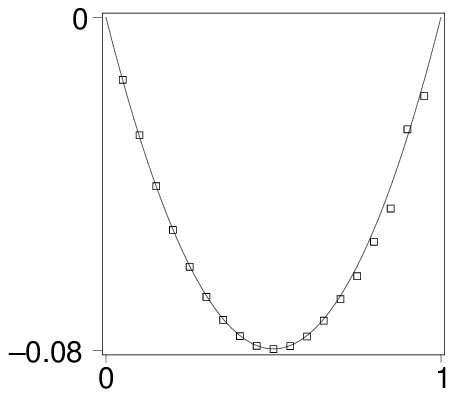}}}
\put(3.43,.9){\footnotesize$\frac{I}{h(I)}$}
\put(4.2,0){\footnotesize$I$}
\end{picture}}}
\caption{\label{101} Numerical verification of the twist predicted
by~\erf{18}. Leftmost: the curved leading edge of the numerically
computed Poincar\'e map indicates a twist map by visibly showing
faster rotation as $I$ increases. The twist overlays a constant
rotation (in $I$) which is caused by high order terms in $a$ and
decreases as $a$ decreases.  Center: the rotation angle per iteration
of the Poincar\'e map on the right as a function of $I$. The slope
corresponds to the twist predicted by~\erf{18}.  Right: reciprocal of
the twist for $I_3=1$, $I_2=I_3-I_1$ as $I_1$ ranges from $0$ to $1$
compared to the parabola predicted by~\erf{18}.}
\end{figure}

The foregoing sort of analysis will always lead to some stability
condition, irrespective of possible errors in the derivation, so it
necessary to check~\erf{18} by comparing it with numerically generated
Poincar\'e maps.  Substitution of $\epsilon=1$ after truncation of
$O(\epsilon)^3$ into~\erf{18} gives the leading behavior of the
Poincar\'e map when $a$, $\nu$, and $I$, are of comparable order (they
are all order $\epsilon^2$). The Poincar\'e map is determined to leading
order in $I$ by the first order twist term $h(I)$ in  when the
zero~order term $2\pi\omega_e/\hat\xi_e^{\alpha_e}$ has a vanishing
effect (i.e. is a multiple of $2\pi$). This happens to occur when
$I_1+I_2=I_3$, as is easily verified. After substitution of
$I_2=I_3-I_1$, the first order twist is
\begin{equation*}
h(I)=-\frac{\pi {I_3}^2}{I_1(I_3-I_1)}I,
\end{equation*}
whereupon
\begin{equation*}
\frac{I}{h(I)}=-\frac1{\pi {I_3}^2}I_1(I_3-I_1),
\end{equation*}
which is a parabola in $I_1/I_3$. As can be seen in Figure~\thrf{101},
this compares well with numerical integrations of the
original (as opposed to  the blown-up)  system.

\section*{Summary}
The following theorem has been proved.
\begin{theorem}
Within the context of the Lagrangian system~\erf{00}, the motion of an
underwater ellipsoid rotating about a long or short principle axis of
inertia is stable modulo $\mtl{SO}(2)\times C$ where $\mtl{SO}(2)$
acts around the rotation axis and $C$ is any cone
containing that axis.
\end{theorem}
This theorem follows from KAM confinement after a blow-up construction
and normal form analysis rather than confinement by Lyapunov functions
derived from energy and momentum.

\section*{Acknowledgments}
This work was supported by an NSERC individual research grant and an
EPSRC visiting fellowship. I thank the University of Warwick
Mathematics Institute for its hospitality during a sabbatical visit
while this paper was written.

\frenchspacing\footnotesize

\end{document}